\documentclass[12pt,leqno]{elsarticle}
\usepackage{amsfonts, color}

\newcounter{conjecture}
\newcounter{remark}

\newtheorem{theorem}{Theorem}
\newtheorem{lemma}{Lemma}

\newcommand{\dd}{\delta}

\newcommand{\lar}{\longrightarrow}
\newcommand{\eps}{\varepsilon}

\newcommand{\reals}{\mathbb{R}}

\newcommand{\nats}{\mathbb{N}}
\newcommand{\lll}{\label}
\newcommand {\rrr}[1]{(\ref{#1})}

\def \be{\begin{equation}}
\def \ee{\end{equation}}
\def \bt{\begin{theorem}}
\def \et{\end{theorem}}
\def \bc{\begin{corollary}}
\def \ec{\end{corollary}}
\def \bl{\begin{lemma}}
\def \el{\end{lemma}}
\def \bea{\begin{eqnarray}}
\def \eea{\end{eqnarray}}
\def \bas{\begin{eqnarray*}}
\def \eas{\end{eqnarray*}}

\def \De{\Delta}

\def \vski{\vspace{12pt}}
\def \ff{\infty}

\def \DD{\Delta}

\def \({\left(}
\def \){\right)}

\def \nn{\nonumber}

\def \bc{\begin{center} }
\def \ec{\end{center} }
\def \bs{\begin{slide} }
\def \es{\end{slide} }

\def\square{{\vcenter{\vbox{\hrule height.3pt
         \hbox{\vrule width.3pt height5pt \kern5pt
            \vrule width.3pt}
         \hrule height.3pt}}}}
\def\qed{{\hfill $\square$ \bigskip}}

\newcounter{cccases}
\begin{document}

\title{Applying Brownian motion to the study of birth-death chains.}

\author{
\begin{tabular}{c}
\textit{Greg Markowsky} \\
gmarkowsky@gmail.com \\
(054)279-5828 \\
Pohang Mathematics Institute \\
POSTECH \\
Pohang, 790-784 \\
Republic of Korea
\end{tabular}}

\begin{abstract}
Basic properties of Brownian motion are used to derive two results concerning birth-death chains. First, the probability of extinction is calculated. Second, sufficient conditions on the transition probabilities of a birth-death chain are given to ensure that the expected value of the chain converges to a limit. The theory of Brownian motion local time figures prominently in the proof of the second result.
\end{abstract}


\maketitle

\section{Introduction}

Let $X_m$ be a Markov chain taking values on the nonnegative integers with the following transition probabilities for $n \neq 0$

\be p_{nj} = \left \{ \begin{array}{ll}
r_{n} & \qquad  \mbox{if } j=n+1  \\
l_n & \qquad \mbox{if } j=n-1 \\ 0 & \qquad \mbox{if } |n-j| \neq 1\;.
\end{array} \right. \ee

Implicit here is the fact that $r_n+l_n=1$. We suppose further for simplicity that $X_0 = k$ almost surely, for some $k \in \nats$. $X_m$ is essentially a random walk on the nonnegative integers, moving to the right from state $n$ with probability $r_n$ and to the left with probability $l_n$. We refer to such a Markov chain as a {\it birth-death chain}. This name comes from considering $X_m$ as the number of members in a population, where at each step either a new member is born or an old member dies, causing the process to increase or decrease by 1. We can assume $p_{00}=1$ and $p_{0j}=0$ for any $j \neq 0$, as when the population reaches 0 it is considered to have gone extinct with no possibility of regeneration. The purpose of this paper is to introduce a method of using properties of Brownian motion to deduce two fundamental theorems concerning birth-death chains. The first theorem, presented in the next section, gives the probability that a birth-death chain goes extinct at some finite time. The second theorem, presented in Section 3, gives sufficient conditions for $E[X_m]$ to converge as $m \lar \ff$. The properties of Brownian motion which will be utilized are standard and can be found in many references on Brownian motion, such as \cite{rosmarc} or \cite{revyor}.

\vski

We will now introduce the basic setup. Let $t_0 := 1$ and

\be \label{pred}
t_n := \frac{l_1 l_2 \ldots l_n}{r_1 r_2 \ldots r_n}
\ee

for $n>0$. Define a sequence $\{x_n\}_{n=0}^\ff$ recursively by setting $x_0=0$, and having defined $x_n$ let $x_{n+1}=x_n+t_n$. Since the sequence $\{x_n\}$ is increasing it converges to a limit $x_\ff$, possibly infinite, as $n \lar \ff$. Let $B_t$ be a Brownian motion starting at $x_k$ and stopped at the first time $T_{\Delta}$ it hits $0$ or $x_\ff$. The recurrence properties of Brownian motion imply that $T_{\De} < \ff$ almost surely. We define a sequence of stopping times $T_m$ which are, roughly speaking, the successive hitting times of $\cal{A} :=$ $ \{ x_n \}_{n=0}^\ff$. More rigorously, $T_m$ is defined recursively by setting $T_0=0$, and having defined $T_m$ we let $T_{m+1}=\inf_{t>T_m}\{B_t \in \cal{A},$ $B_t \neq B_{T_m}\}$. We see that the variables $B_{T_0},B_{T_1},B_{T_2}, \ldots $ form a random process taking values in $\cal{A}$. The strong Markov property of Brownian motion, together with the standard exit distribution of Brownian motion from an interval, imply that

\bea \label{}
\nn && P(B_{T_{m+1}} = x_{n+1} | B_{T_{m}} = x_{n}) = \frac{x_{n}-x_{n-1}}{x_{n+1}-x_{n-1}} = \frac{t_{n-1}}{t_{n-1}+t_n} = \frac{1}{1+l_n/r_n} = r_n
\eea

and, likewise,

\bea \label{}
\nn && P(B_{T_{m+1}} = x_{n-1} | B_{T_{m}} = x_{n}) = l_n
\eea

If we define $\phi$ on $\{x_n\}_{n=0}^\ff$ by $\phi(x_n)=n$, we see that $\phi(B_{T_0}),\phi(B_{T_1}),\phi(B_{T_2}), \ldots $ is a realization of our original birth-death chain. The picture below gives an example, where we have oriented the time axis vertically and the space axis horizontally.

\vski

\hspace{.8cm} \includegraphics[width=110mm,height=80mm]{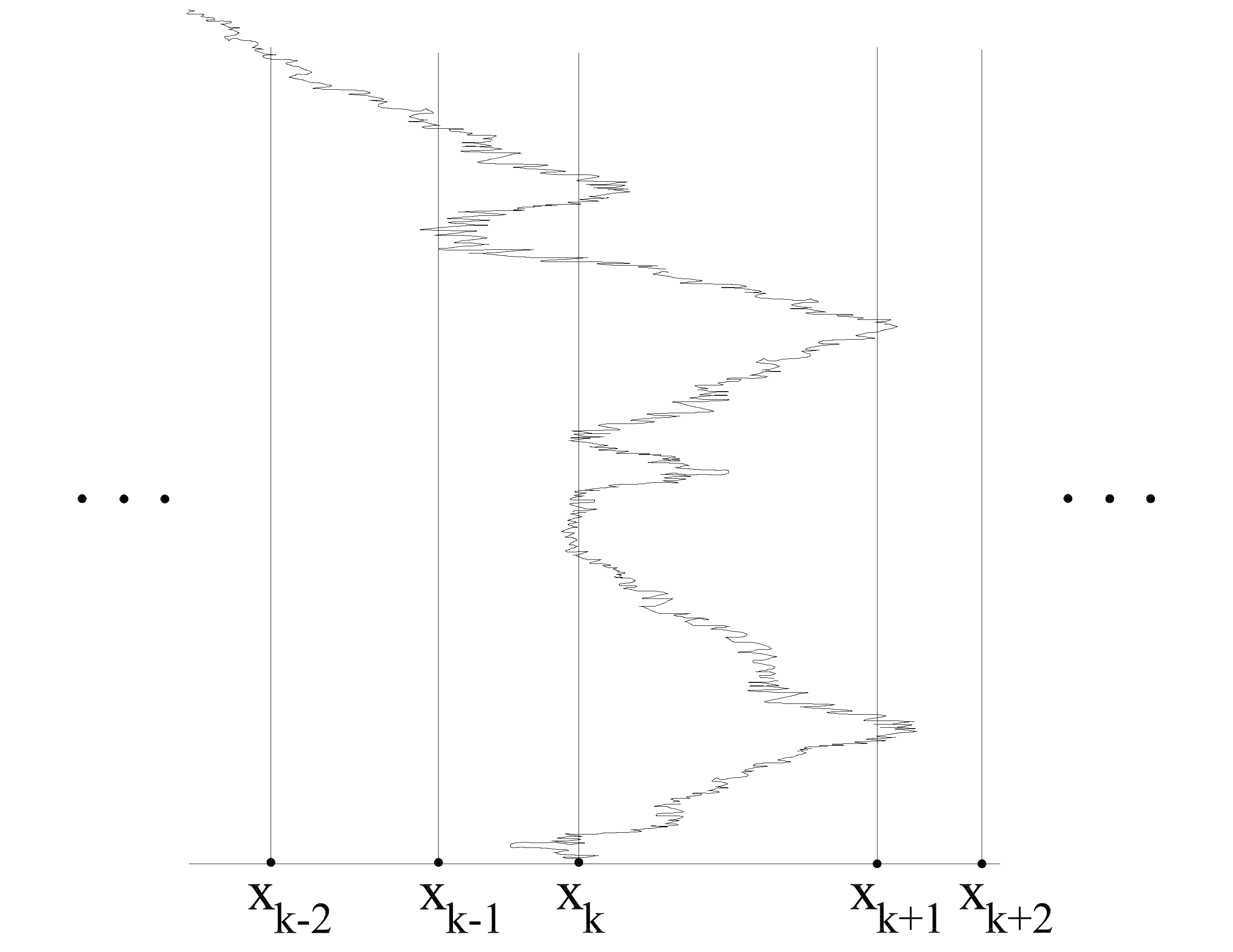}

{\small Figure 1: The Brownian path pictured realizes the birth-death path $k,k+1,k,k+1,k,k-1,k,k-1,k-2, \ldots$}

\vski

Given this framework, we are ready to prove several theorems. In the sequel, any reference to $X, B, x_n, T_{\Delta},\phi,$ etc. will refer to the definitions presented in this section.

\section{The extinction probability of a birth-death chain}

Perhaps the most fundamental question one can ask regarding a birth-death chain is whether the population must go extinct or not, that is, whether $P(X_m = 0$ for some $m)=1$ or $P(\lim_{m \lar \ff} X_m = +\ff)>0$. Let $P_k$ be the probability that the birth-death chain eventually hits 0 (recall $X_0=k$ a.s.). We then have the following.

\bt \lll{surf}
\be \label{mass}
P_k = \frac{\sum_{j=k}^\ff t_j}{\sum_{j=0}^\ff t_j}
\ee

where this quotient is interpreted as being equal to 1 if the sums diverge.
\et

This elegant theorem has a straightforward proof using recurrence relations; see \cite{norr} or \cite{sysk}. A potentially pleasing aspect of the proof below, however, lies in giving a clear, visual intuition for the sums in \rrr{mass}.

\vski

{\bf Proof of Theorem \ref{surf}:}
Recall that $x_\ff = \lim_{n \lar \ff} x_n$ is given by

\be \label{}
x_\ff=\sum_{j=0}^\ff t_j
\ee

If $x_\ff=\ff$, so that both sums in \rrr{mass} diverge, then $B_{T_{\DD}}=0$ almost surely. This implies that the population dies out with probability $1$. On the other hand, if $x_\ff<\ff$ then $P(B_{T_{\De}}=0)$ is given by

\be \label{}
\frac{x_\ff-x_k}{x_\ff-0} = \frac{\sum_{j=k}^\ff t_j}{\sum_{j=0}^\ff t_j}
\ee

However, as in the first case, $P(B_{T_\De}=0)$ is precisely $P_k$, the probability of extinction. This is because the Brownian motion hitting $x_\ff$ before $0$ implies $B_{T_m}\lar x_\ff$, hence $\phi(B_{T_m}) \lar \ff$, whereas hitting $0$ before $x_\ff$ implies $\phi(B_{T_\De}) =0$ for some $m$. The two cases ($x_\ff=\ff$ and $x_\ff < \ff$) are illustrated in the following figure.

\vspace{-1.5cm}

\includegraphics[width=140mm,height=110mm]{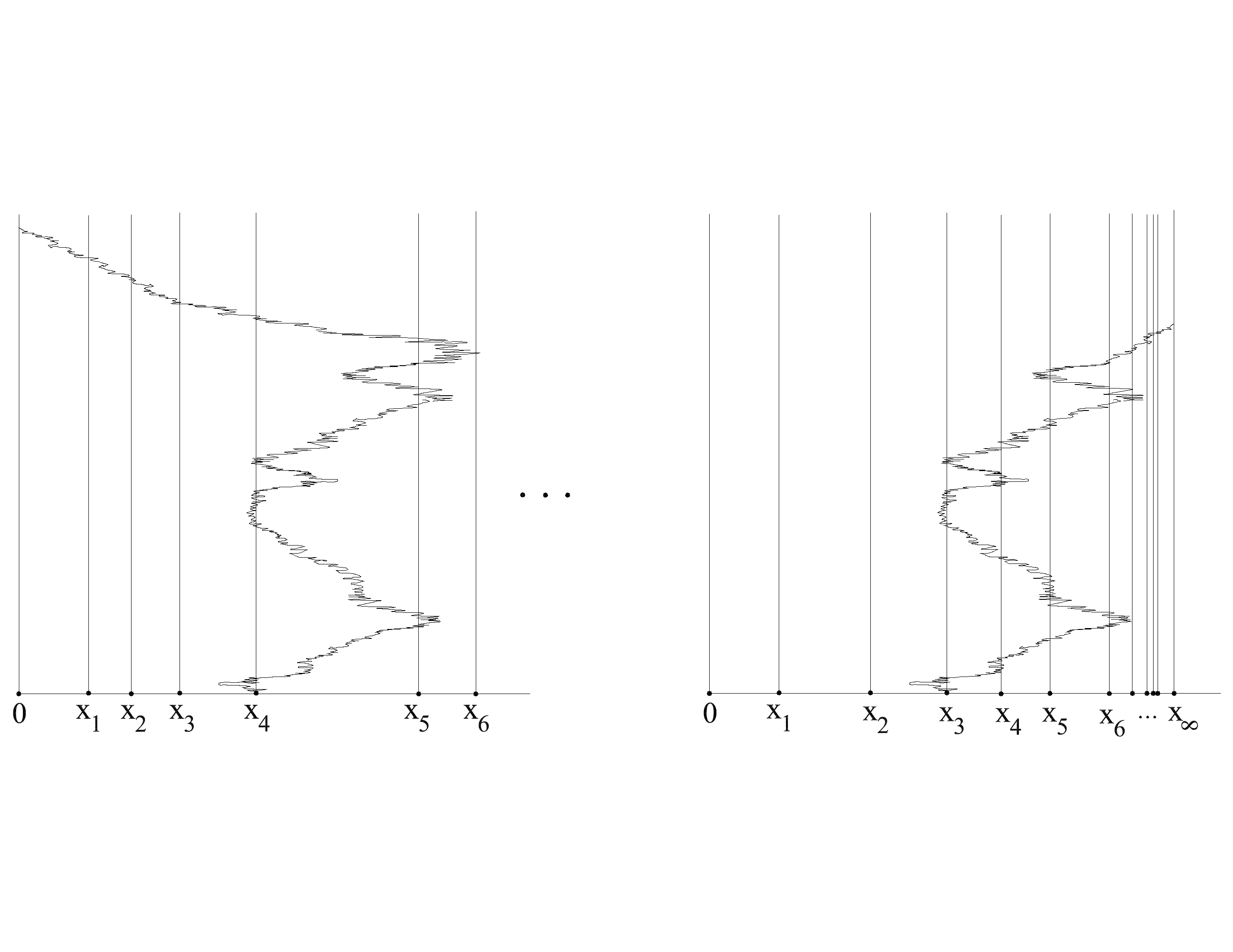}

\vspace{-2cm}
{\small Figure 2: The left panel illustrates the situation in which $\sum_{j=0}^\ff t_j$ diverges. Eventually, $B_t$ hits $0$ and the population goes extinct. The right panel illustrates the other scenario, in which $\sum_{j=0}^\ff t_j = x_\ff < \ff$. In this case, there is a positive probability that $B_t$ hits $x_\ff$ before $0$, in which case the population never goes extinct.}

\vski

This completes the proof of Theorem 1. \qed

\section{The long-term average of a birth-death chain}

Recall that $B$ and $X$ are stopped upon reaching $0$. It will therefore be convenient to let $X_m$ be defined to be $0$ for all $m>m_0$, where $m_0$ is the smallest integer, if it exists, for which $X_{m_0}=0$. Similarly, for convenience let $T_m=T_\DD$ for all $m > m_0$, where $m_0$ is the smallest integer, if it exists, for which $B_{T_{m_0}}=0$. In the case $r_i=l_i=\frac{1}{2}$ for all $i$, it is well known that $X_m$ is a martingale, and therefore $E[X_m]=E_0=k$ for all $m$. This occurs despite the fact that $P(X_m = 0) \lar 1$ as $m \lar \ff$, as the average value of $X_m$ on $\{X_m \neq 0\}$ grows at exactly the right speed to balance the set of large probability upon which $X_m=0$. Such behavior certainly does not hold for the general case, since we no longer have the martingale property, but we will see that the Brownian motion model presented above can shed light on the behavior of $E[X_m]$ as $m \lar \ff$.

\vski

Recall that $\cal{A}$ $=\{x_n\}_{n=0}^\ff$. Let $\phi: \cal{A}$ $\lar \reals^+$ be extended to a continuous function from $\reals^+$ to $\reals^+$ by defining $\phi$ to be linear on each interval $(x_{n-1},x_n)$. Alternatively, we may think of $x_n=x(n)$ as a function from $\mathbb{N}$ to $\reals$ which can be extended by linear interpolation to an increasing function from $\reals^+$ to $\reals^+$. In this case, $\phi$ is simply $x^{-1}$. $\phi$ is therefore a piecewise linear function, and $\phi'$ exists on $\reals^+ - \cal{A}$. Let $\phi'_n$ be the value of $\phi'$ on $(x_{n-1},x_n)$. We will prove the following theorem.

\bt \lll{bigguyii}
If $\phi'_\ff = \lim_{n \lar \ff} \phi'_n$ exists, then

\be \lll{yokoneg1}
\lim_{m \lar \ff} E[X_m] = x_k \phi'_\ff
\ee

\et

Note that we are allowing $\phi'_\ff = +\ff$ or $0$. Writing $\phi'_\ff$ and $x_k$ in terms of the $l_n$'s and $r_n$'s shows that the following statement is equivalent to Theorem \ref{bigguyii}.

\vski

{\it If $t_\ff := \lim_{n \lar \ff} \frac{l_1\ldots l_{n}}{r_1\ldots r_{n}}$ exists then $\lim_{m \lar \ff}E[X_m]$ exists, and }

\be \lll{}
\lim_{m \lar \ff}E[X_m] = \frac{1+ \frac{l_1}{r_1} + \ldots + \frac{l_1\ldots l_{k-1}}{r_1\ldots r_{k-1}}}{t_\ff}
\ee

\vski

The bulk of the rest of this section is devoted to the proof of this theorem. We will simplify initially by assuming $\sum_{n=1}^{\ff} |\phi'_{n+1}-\phi'_n| < \ff$; this condition will be removed at the end of the proof. For the case in which there is a positive probability that the population never goes extinct, it is easy to see that $E[X_m] \lar \ff$ as $m \lar \ff$, and that $\phi'_\ff$ exists and is equal to $+\ff$, so that \rrr{yokoneg1} is valid. We will therefore assume that $P(X_m = 0$ for some $m) = 1$. Note that $\phi'_{n+1} = \frac{1}{t_n}$, and $x_{n+1}-x_n = t_n$, so that $\phi_{n+1}'(x_{n+1}-x_n) = 1$. Note also that $x_1 \phi'_1 = 1$. This allows us to perform the following manipulations to obtain an expression which will be more convenient for the purposes of the proof.

\bea \lll{yokoa}
&& x_k \phi'_\ff = x_k (\phi'_\ff - \phi'_{k}) + x_k \phi'_k + k - x_1\phi'_1 - \sum_{n=1}^{k-1} \phi_{n+1}(x_{n+1}-x_n)
\\ \nn && \hspace{1.1cm}= k + x_k (\phi'_\ff - \phi'_{k}) + \sum_{n=1}^{k-1} (\phi'_{n+1}-\phi'_n) x_n
\eea

The last equality uses summation by parts; see \cite{lang}. We see that the conclusion of the theorem is equivalent to showing

\be \lll{yoko}
\lim_{m \lar \ff} E[X_m] = k + x_k (\phi'_\ff - \phi'_{k}) + \sum_{n=1}^{k-1} (\phi'_{n+1}-\phi'_n) x_n
\ee

This is what we will prove. We will proceed through several lemmas, and will need properties of Brownian motion {\it local time}, which is the density of the occupation measure of Brownian motion with respect to Lebesgue measure. That is, the local time $L_t^x$ satisfies

\be \lll{}
L_t^x dx = \int_{0}^{t} 1_{B_s \in dx} ds
\ee

It is well known that $L_t^x$ exists and that

\be \lll{}
L_t^x = \lim_{\eps \lar 0} \frac{1}{2\eps} \int_{0}^{t} 1_{|B_s-x|<\eps}ds
\ee

almost surely. See \cite{rosmarc} for a comprehensive treatment of local time, or the more general reference \cite{revyor}. The following is an extension of Tanaka's formula, Theorem VI.1.2 in \cite{revyor}.

\bl \lll{45}
Almost surely, for any stopping time $T$,

\be \lll{maz}
\phi(B_T) = k + \int_{0}^{T} \phi'(B_s) dB_s + \sum_{n=1}^\ff \frac{(\phi'_{n+1}-\phi'_n)}{2} L^{x_n}_T
\ee

where $L^{x_n}_T$ denotes the local time of $B_t$ at $x_n$ at time $T$.
\el

{\bf Proof:} Note that $\phi''(x) = \sum_{n=1}^\ff (\phi'_{n+1}-\phi'_n) \dd_{x_n}(x)$ in the sense of distributions, where $\dd_{x_n}(x) = \dd_0(x-x_n)$ denotes the Dirac delta function at point $x_n$. Lemma \ref{45} is therefore seen to be a special case of the It\^{o}-Tanaka formula, Theorem VI.1.5 of \cite{revyor}, provided that $\phi$ can be realized as a difference of two convex functions. However, any piecewise-linear function can be realized as the difference of two convex functions, provided that the points of nondifferentiability do not accumulate. We may argue as follows. $\phi'$ is a piecewise constant function, which is therefore of bounded variation on bounded intervals, and as such we may write $\phi' = f-g$ where $f$ and $g$ are nondecreasing. Let $F$ and $G$ be antiderivatives of $f$ and $g$ chosen so that $\phi = F-G$. Then $F$ and $G$ are convex, and the result follows. \qed

Applying this lemma to the stopping time $T_m$, we immediately obtain

\be \lll{kok}
E[X_m] = k + \sum_{n=1}^\ff \frac{(\phi'_{n+1}-\phi'_n)}{2} E[L^{x_n}_{T_m}],
\ee

Note that the convergence of the sum at this point is not an issue, since $L^{x_n}_{T_m} = 0$ for $n>m+k$. Using the identity \rrr{kok} does not seem to be an effective way to calculate $E[X_m]$, due to the difficulty of obtaining information about $T_m$. Nonetheless, we do know that $T_m \nearrow T_\DD$ as $m \lar \ff$, and this implies

\be \lll{oko}
\lim_{m \lar \ff} E[X_m] = k + \sum_{n=1}^\ff \frac{(\phi'_{n+1}-\phi'_n)}{2} E[L^{x_n}_\ff],
\ee

provided that $E[L^{x_n}_\ff]$ can be bounded uniformly, which we will show soon to be the case. We should mention that it was in obtaining \rrr{oko} that we used the assumption that $\sum_{n=1}^{\ff} |\phi'_{n+1}-\phi'_n| < \ff$. This is because a priori the quantities $E[L^{x_n}_{T_m}]$ may be growing in some strange way that causes problems if $\sum_{n=1}^{\ff} |\phi'_{n+1}-\phi'_n| = \ff$. We will return to this point at the end of the proof. In light of \rrr{oko}, we must compute $E[L^{x_n}_\ff]$.

\bl \lll{sbtm}
\be \lll{a3}
E[L^{x_n}_\ff] = 2 \min (x_k,x_n)
\ee
\el

{\bf Proof:} One may derive this through standard calculations involving the probability density function of $B_t$, but the following is a quicker and easier proof. Let us suppose first that $n=k$. From Tanaka's formula, $E[L^{x_k}_\ff] = \lim_{t \lar \ff} E[|B_t-x_k|]$. Furthermore, $B_t$ is a martingale, so $E[(B_t-x_k)]=0$ for all $t$. It follows from this that

\bea \label{}
&& E[L^{x_k}_\ff]=2\lim_{t \lar \ff} E[\max(-(B_t-x_k),0)]
\\ \nn && \hspace{1.3cm} = 2 \lim_{t \lar \ff} \Big( x_k P(B_t=0) + \int_{0}^{k} (x_k-x) P(B_t \in dx) \Big)
\eea

As $\lim_{t \lar \ff} P(B_t=0) = 1$, we can conclude that $E[L^{x_k}_\ff]=2x_k$. Now suppose $n \neq k$. Note that, if we let $T_{x_n} = \inf \{t: B_t=x_n\}$, then

\be \lll{runy}
L^{x_n}_\ff = L_{T_{x_n}} + L_\ff(B \circ \theta_{T_{x_n}})1_{T_{x_n}<T_\DD} = L_\ff(B \circ \theta_{T_{x_n}})1_{T_{x_n}<T_\DD}
\ee

where $\theta$ denotes the standard shift operator and $L_t(B \circ \theta_{T_{x_n}})$ is the local time of the shifted process $B \circ \theta_{T_{x_n}}$. 
Let $E_{x_j}$ denote expectation with respect to a Brownian motion $W$ which starts at $x_j$ and is stopped upon hitting $0$. The prior calculation together with \rrr{runy} and the strong Markov property of Brownian motion imply that

\bea \label{}
&& E[L^{x_n}_\ff] = P(T_{x_n}<T_\DD) E_{x_n} [L_\ff^{x_n}]
\\ \nn && \hspace{1.32cm} = P(T_{x_n}<T_\DD) 2x_n
\eea

The general result follows from noting that $P(T_{x_n}<T_\DD)$ is $1$ if $x_n < x_k$ and $\frac{x_k}{x_n}$ if $x_n > x_k$. \qed

Combining \rrr{oko} and Lemma \ref{sbtm} gives

\be \lll{yoko2}
\lim_{m \lar \ff} E[X_m] =  k + \sum_{n=1}^{k-1} (\phi'_{n+1}-\phi'_n) x_n + x_k \sum_{n=k}^{\ff}(\phi'_{n+1} - \phi'_n)
\ee

Since $\sum_{n=k}^{\ff}(\phi'_{n+1} - \phi'_n) = (\phi'_\ff-\phi'_{k})$, we are done in this case. It remains only to remove the restriction that $\sum_{n=1}^{\ff} |\phi'_{n+1}-\phi'_n| < \ff$. The following lemma is key.

\bl \lll{}
For any $m$, and any $n \geq k$, $E[L^{x_n}_{T_m}] \geq E[L^{x_{n+1}}_{T_m}]$.
\el

{\bf Proof:} In fact, we may prove somewhat more, namely that if $T$ is any stopping time with $B_T \in \cal{A}$ almost surely, then $E[L^{x_n}_{T}] \geq E[L^{x_{n+1}}_{T}]$. Let $E_{x_j}$ and $W$ be as in the proof of Lemma \ref{sbtm}.
Using Lemma \ref{sbtm} and the strong Markov property of Brownian motion, we obtain

\bea \label{}
&& E[L^{x_n}_{T}] = E[L^{x_n}_{\ff}] - E\Big[\lim_{\eps \lar 0}\frac{1}{2\eps}\int_{T}^{\ff} 1_{(-\eps,\eps)}(B_s - x_n) ds \Big]
\\ \nn && \hspace{1.33cm} = 2 x_k - \sum_{j=1}^{\ff} P(B_T = x_j)E_{x_j} \Big[ \lim_{\eps \lar 0} \frac{1}{2\eps} \int_{0}^{\ff} 1_{(-\eps,\eps)}(W_s- x_n) ds \Big]
\\ \nn && \hspace{1.33cm} = 2 x_k - \sum_{j=1}^{n-1} P(B_T = x_j)x_j - \sum_{j=n}^{\ff} P(B_T = x_j)x_n
\eea

Similarly,

\be \lll{}
E[L^{x_{n+1}}_{T}] = 2 x_k - \sum_{j=1}^{n} P(B_T = x_j)x_j - \sum_{j=n+1}^{\ff} P(B_T = x_j)x_{n+1}
\ee

The conclusion of the lemma now follows from the fact that $x_{n+1} > x_n$. \qed

We may now complete the proof of the theorem. Recall \rrr{kok}, and observe that $E[L^{x_n}_{T_m}] = 0$ for $n>k+m$, since $X_m \leq m+k$. This means that \rrr{kok} is in fact a finite sum.

\bea \label{kok2}
&& E[X_m] = k + \sum_{n=1}^{m+k} \frac{(\phi'_{n+1}-\phi'_n)}{2} E[L^{x_n}_{T_m}]
\\ \nn && \hspace{1.3cm} = k + \sum_{n=1}^{k-1} \frac{(\phi'_{n+1}-\phi'_n)}{2} E[L^{x_n}_{T_m}] + \sum_{n=k}^{k+m} \frac{(\phi'_{n+1}-\phi'_n)}{2} E[L^{x_n}_{T_m}]
\eea

The indices of the first sum in the final expression of \rrr{kok2} are independent of $m$. This implies that the sum converges as $m\lar \ff$, since $E[L^{x_n}_{T_m}] \lar 2x_n$ as $m \lar \ff$ for $n \leq k$. We must show that the second sum converges as $m\lar \ff$. We use summation by parts again, which gives

\bea \label{kok2}
&& \sum_{n=k}^{k+m} \frac{(\phi'_{n+1}-\phi'_n)}{2} E[L^{x_n}_{T_m}] = \frac{1}{2}\Big( \phi'_{k+m+1}E[L^{x_{k+m+1}}_{T_m}] - \phi'_{k}E[L^{x_{k}}_{T_m}]
\\ \nn && \hspace{1cm}  - \sum_{n=k}^{k+m} \phi'_{n+1} (E[L^{x_{n+1}}_{T_m}] - E[L^{x_{n}}_{T_m}]) \Big)
\eea

Let us assume that $\phi'_\ff < \ff$, and let $\eps>0$ be given. We may choose $N>k$ such that $\phi'_n \in (\phi_\ff-\eps,\phi_\ff+\eps)$ for all $n \geq N$. Having chosen this, we may choose $M > N-k$ such that $2x_k \geq E[L^{x_{n}}_{T_m}] > 2x_k - \eps$ for all $n \in [k,N], m \geq M$. Using the fact that $E[L^{x_{k+m+1}}_{T_m}]=0$, and setting $\overline{\phi'} = \sup_{j > 0} \phi'_j$, we see that for $m>M$

\bea \label{}
&& \sum_{n=k+1}^{k+m} \frac{(\phi'_{n+1}-\phi'_n)}{2} E[L^{x_n}_{T_m}]
\\ \nn && \hspace{1cm} \leq \frac{1}{2}\Big( -\phi'_{k} (2x_k - \eps) + \sum_{n=k}^{N} \phi'_{n+1} (E[L^{x_{n}}_{T_m}] - E[L^{x_{n+1}}_{T_m}])
\\ \nn && \hspace{2cm} + \sum_{n=N+1}^{k+m} \phi'_{n+1} (E[L^{x_{n}}_{T_m}] - E[L^{x_{n+1}}_{T_m}]) \Big)
\\ \nn && \hspace{1cm} \leq \frac{1}{2}\Big( -\phi'_{k} (2x_k - \eps) + \overline{\phi'} (E[L^{x_{k}}_{T_m}] - E[L^{x_{N+1}}_{T_m}])
\\ \nn && \hspace{2cm} + (\phi'_\ff + \eps) (E[L^{x_{N+1}}_{T_m}] - E[L^{x_{k+m+1}}_{T_m}]) \Big)
\\ \nn && \hspace{1cm} \leq \frac{1}{2}\Big( -\phi'_{k} (2x_k - \eps) + \overline{\phi'} \eps + (\phi'_\ff + \eps) 2x_k \Big)
\eea

This shows that

\be \lll{}
\limsup_{m \lar \ff} \sum_{n=k}^{k+m} \frac{(\phi'_{n+1}-\phi'_n)}{2} E[L^{x_n}_{T_m}] \leq x_k(\phi'_\ff - \phi'_{k})
\ee

Proceeding similarly, we can obtain

\be \lll{}
\liminf_{m \lar \ff} \sum_{n=k+1}^{k+m} \frac{(\phi'_{n+1}-\phi'_n)}{2} E[L^{x_n}_{T_m}] \geq x_k(\phi'_\ff - \phi'_{k})
\ee

Together, these prove the desired convergence. The case $\phi'_\ff = +\ff$ is similar but easier and is omitted. This completes the proof of Theorem \ref{bigguyii}. \qed

We conclude with a simple but counterintuitive example. Let $l_n=\frac{n}{2n+1}, r_n = \frac{n+1}{2n+1}$ for $n \geq 1$. Then $t_n = \frac{1}{n+1}$, so that $t_\ff = 0$. On the other hand, $x_\ff = 1 + \sum_{n=1}^{\ff}t_n = \ff$. We see that the birth-death chain $X_m$ built upon these transition probabilities has an extinction probability of 1, but $E[X_m] \lar \ff$ as $m \lar \ff$.

\section{Acknowledgements}

I am grateful for support from the Priority Research Centers Program through the National Research Foundation of Korea (NRF) funded by the Ministry of Education, Science and Technology (Grant \#2009-0094070). I would also like to thank the referee for comments which improved the exposition, as well as for asking an interesting question which led to the development of Section 3.

\end{document}